\documentclass[12pt, reqno]{amsart}

\usepackage{amsfonts}
\usepackage{amsmath}
\usepackage{amsthm}
\usepackage{amssymb}
\usepackage{latexsym}
\usepackage{enumerate}
\usepackage{fixmath}
\usepackage{upgreek}
\usepackage{hyperref}
\usepackage{tabularx}
\usepackage{amscd,amssymb}
\input xy
\xyoption{all}

\advance\textwidth by 1.2in \advance\oddsidemargin by -.6in
\advance\evensidemargin by -.6in

\newtheorem{lem}{Lemma}
\newtheorem{prop}{Proposition}
\newtheorem*{defn}{Definition}
\newtheorem{thm}{Theorem}
\newtheorem*{theorem}{Theorem}

\newcommand{\lie}[1]{\mathfrak{#1}}
\newcommand{\mat}[1]{\left[ \begin{matrix} #1 \end{matrix} \right]}

\newenvironment{pf}{\proof}{\endproof}
 \newcommand\bpi{{\mbox{\boldmath $\pi$}}}

\newcommand{\twistg}{L(\lie g, \sigma, m)}
\newcommand{\twisth}{L(\lie h, \sigma, m)}
\renewcommand{\vec}{\mathbf}
\newcommand{\aut}{\operatorname{aut}}
\newcommand{\Id}{\operatorname{Id}}
\newcommand{\ad}{\operatorname{ad}}
\newcommand{\ev}{\operatorname{ev}}
\newcommand{\supp}{\operatorname{supp}}
\newcommand{\nilrad}{\operatorname{nilrad}}
\usepackage{hyperref} 
\title[Finite--dimensional representations of loop algebras]{Finite--dimensional representation theory of loop algebras: a survey}
\author{Prasad Senesi}
\address{Prasad Senesi, Department of Mathematics and Statistics, University of Ottawa, Ottawa, ON, Canada, K1N 6N5.}

\begin{document}

\maketitle
\thispagestyle{empty}

\begin{abstract}

We survey some important results concerning the finite--dimensional representations of the loop algebras $\lie g \otimes \mathbb{C}\left[ t^{\pm 1} \right]$ of a simple complex Lie algebra $\lie g$, and their twisted loop subalgebras.  In particular, we review the parametrization and description of the Weyl modules and of the irreducible finite--dimensional representations of such algebras, describe a block decomposition of the (non--semisimple) category of their finite--dimensional representations, and conclude with recent developments in the representation theory of multiloop algebras.
\end{abstract}

\tableofcontents

\section{Introduction}

In this survey we review some important developments concerning the theory of finite--dimensional representations of the \textit{loop algebras}, a class of infinite--dimensional Lie algebras.  These Lie algebras come in one of two varieties, either \textit{untwisted} or \textit{twisted}.  The untwisted Lie algebras are all of the form $L(\lie g) = \lie g \otimes \mathbb{C}\left[ t ^\pm \right]$, where $\lie g$ is a simple finite--dimensional Lie algebra.  The twisted loop algebras all occur as certain subalgebras of the untwisted algebras; their precise description is found in Section \ref{Loop algebras}.  

Loop algebras occur in a realization of the affine Kac-Moody algebras given by V. Kac in \cite{Kac}.  The representation theory of these affine algebras has gathered much interest since their introduction 40 years ago, and is related to a variety of algebraic and geometric topics in mathematics and in mathematical physics, including crystal bases, vertex operator algebras, conformal field theory, solvable lattice models, and solutions to the Yang--Baxter equation.  The loop algebras also provide important examples of \textit{centreless Lie tori}, which in turn are essential in the description of the \textit{extended affine Lie algebras}, or EALAs.  These are generalizations of both the finite--dimensional and affine Kac--Moody algebras, and their structure and classification is the subject of many recent papers; see, for example, \cite{AABGP, ABFP1, ABFP2, N1, N2}. 

In 1986, V. Chari classified certain irreducible infinite--dimensional representations of the \textit{extended} (untwisted) loop algebra $L(\lie g)^e$ (see Section \ref{additional}), and with A. Pressley these representations were described and the classification was extended to the twisted case; see \cite{C, CPunitary, CPint}.  In these references, the (tensor products of) evaluation representations of $L(\lie g)$ are defined, and then used in the construction of infinite--dimensional $L(\lie g)^e$--modules.  It is shown in \cite{CFS, CPweyl, R} that these evaluation modules are precisely the irreducible representations of the loop algebras.

As we should expect from its description, the representation theory of $\lie g$ is very helpful in understanding that of the corresponding loop algebra.  Indeed, the concepts of highest--weight and Verma modules of $\lie g$ play a prominent role in the results discussed here, and the dominant integral linear forms of a Cartan subalgebra $\lie h$ of $\lie g$ occur in the parametrization of the simple finite--dimensional representations of the loop algebra, among other things.  

There are four themes which we emphasize here.  The first is the structure of the twisted loop algebras.  This structure is dependent upon the choice of a finite--order Lie algebra automorphism of $\lie g$, and this choice in turn determines a decomposition of the loop algebra with respect to which the highest--weight theory of $\lie g$ can be extended to the loop case.  A theorem of V. Kac guarantees that, up to isomorphism, it is sufficient to consider only the twisted loop algebras constructed from \textit{diagram automorphisms} of $\lie g$.  This restriction cannot be made, however, when we consider representations of the multiloop algebras (see Section \ref{multiloops}).  The second theme is the classification of the finite--dimensional irreducible representations of a twisted loop algebra.  This classification consists of two problems: describing all finite--dimensional irreducible representations, and providing a parametrization of their isomorphism classes.  An approach to the first problem involving a careful analysis of the kernels of irreducible representations (see \cite{B, Lau, R}) is reviewed in Section \ref{multiloops}.  The third theme is the concept of a Weyl module of a loop algebra, first defined in \cite{CPweyl} and extended to the twisted case in \cite{CFS}.  This is a loop algebra analogue of the maximal highest--weight Verma module used in the representation theory of $\lie g$ (although the Weyl modules are always finite--dimensional). Just as occurs in the case of the Verma modules, the Weyl modules are in bijective correspondence with the (isomorphism classes of) finite--dimensional irreducible modules -- this correspondence given by taking quotients of Weyl modules by their unique maximal ideals.  The parametrization of the Weyl modules is achieved with a multiplicative monoid $\mathcal{P}$ of $\text{rank}(\lie g)$--tuples of polynomials (sometimes called the \textit{Drinfeld polynomials} in the literature).  The fourth theme is an examination of the category of all finite--dimensional representations of a loop algebra.  This category contains representations of the loop algebra which are indecomposable, yet reducible.  Therefore not every object can be written as a sum of irreducible objects.  However, the category itself will decompose into a sum of indecomposable categories.  Such a decomposition has the property that each indecomposable object will
lie in a unique indecomposable subcategory, although such a subcategory may contain many
nonisomorphic indecomposables. In such a case, when complete reducibility is not at hand, it is natural to
search for a description of a maximal decomposition of the category called a \textit{block decomposition}. This is a familiar and useful strategy in the
BGG category $\mathcal{O}$, for example, where the blocks are parametrized by central characters of the universal
enveloping algebra of $\lie g$.  It was proved in \cite{CM} and \cite{S} that the blocks of the category of finite--dimensional representations of a twisted loop algebra are parametrized by certain equivalence classes of elements of $\mathcal{P}$ called the (twisted) \textit{spectral characters}.

Perhaps the most natural direction to proceed after these issues have been resolved is to generalize the loop algebra $\lie g \otimes \mathbb{C}\left[ t^\pm \right]$ by replacing $\mathbb{C}\left[ t^\pm \right]$ with the ring of Laurent polynomials in more than one variable $\mathbb{C}\left[ t_1^\pm, \ldots, t_N^\pm \right]$.  The corresponding Lie algebras $\lie g \otimes \mathbb{C}\left[ t_1^\pm, \ldots, t_N^\pm \right]$ and their (suitably defined) twisted subalgebras are known as \textit{multiloop algebras}.  We continue with a section describing recent results (\cite{Lau}) by M. Lau in which the irreducible modules of these multiloop algebras are classified.  We then conclude with a small -- and by no means exhaustive -- selection of research pursuits closely related to the finite--dimensional representation theory described here. 

\noindent{\bf Acknowledgments:}  I extend all of my gratitude to the following persons: to V. Chari, who was and remains my teacher, for her continuing guidance and support; to G. Fourier, with whom V. Chari and myself collaborated in the classification of the twisted Weyl modules; to A. Moura, for his support and help in understanding the block decompositions described here; to M. Lau, for his very important recent contributions and discussions concerning the multiloop algebras, and to E. Neher and A. Savage, for their tireless proofreading and helpful suggestions, which are too numerous to list.

\section{Loop algebras}\label{Loop algebras}

In this section we will describe the Lie algebras whose representations are the subject of the paper.  These Lie algebras are all referred to as loop algebras, or sometimes `loops in $\frak g$'.  Their universal central extensions are the derived algebras of the affine Kac--Moody algebras.    

\pagebreak 

\subsection{Preliminaries}\mbox{}

All vector spaces and representations are defined over the ground field $\mathbb{C}$.  If $\lie a$ is a Lie algebra, $U(\lie a)$ is its universal enveloping algebra.  We denote by $\lie g$ a simple finite--dimensional Lie algebra, and $\lie h$ a Cartan subalgebra of $\frak g$.   For $\alpha \in \frak h^*$, let $\frak g_\alpha = \left\{ x \in \frak g \left. \right| \left[ h ,x \right] = \alpha(h) x, \; h \in \frak h \right\}$. 

If $\lie g$ is a simple finite--dimensional Lie algebra with Cartan matrix $X = (X_{ij})$ (i.e., $\lie g$ is of type $X$), we will write $\lie g = \lie g(X)$.  Let $R$ be the set of weights of the pair ($\frak g$, $\lie h$), i.e., $ R = \left\{ \alpha \in \frak h^* \left. \right| \alpha \neq 0; \; \;  \frak g_\alpha \neq 0 \right\}$. 
 We set $n = \text{rank}(\lie g)$, choose a set of simple roots $\Delta = \left\{ \alpha_1, \ldots, \alpha_n \right\} \subseteq \lie h^*$, let $Q^+ = \bigoplus_{i=1}^n \mathbb{Z}_{\geq 0}\alpha_i$, $Q^- = - Q^+$,  $Q = Q^-  + Q^+$, 
and $R^\pm = R \cap Q^\pm$.  We then have a triangular decomposition of the triple $(\frak g, \lie h, \Delta)$
 \[ \frak g = \frak n^- \oplus \frak h \oplus \frak n^+,\]
 where $\frak n^\pm = \bigoplus_{\alpha \in R^\pm}\frak g_\alpha$.  We remind the reader here that this decomposition of $\frak g$ is dependent upon our choice of Cartan subalgebra $\frak h \subseteq \frak g$ and of simple roots $\Delta$.  For example, when $\frak g \cong \frak{sl}_2$, the standard choice of Cartan subalgebra is
  \[\lie h = \mathbb{C}\mat{1&0\\0&-1}, \text{for which we have } \lie n^+ = \mathbb{C}\mat{0&1\\0&0}, \; \lie n^- = \mathbb{C}\mat{0&0\\1&0}.\]
However, we may also take for a Cartan subalgebra 
\[ \lie h = \mathbb{C}\mat{0&1\\-1&0},\]
in which case the triangular decomposition of $\frak g$ is 
\[ \frak g = \mathbb{C}\mat{i&1\\1&-i} \oplus \mathbb{C}\mat{0&1\\-1&0} \oplus \mathbb{C}\mat{-i&1\\1&i}.\]

 Let $\left\{ h_1, \ldots, h_n \right\} \subseteq \frak h$ be the vector space basis of $\lie h$ such that $\alpha_i(h_j) = X_{ij}$ (a set of coroots), and denote by $\omega_1, \ldots, \omega_\ell$ the corresponding set of fundamental weights in $\frak h^*$; i.e., those linear forms satisfying $\omega_i(h_j) = \delta_{ij}$; set $P^+ = \bigoplus_{i=1}^n \mathbb{Z}_{\geq 0}\omega_i$, $P^- = - P^+$,  $P = P^-  + P^+$.

Let $\text{Aut}(\lie g)$ be the group of Lie algebra automorphisms of $\lie g$, and $\text{aut}(\lie g)$ the subgroup of finite--order Lie algebra automorphisms.  In the next section we will construct Lie algebras from $\lie g$ using elements of $\text{aut}(\lie g)$.  

\subsection{The loop algebras}\label{loopdef}\mbox{}

Let $\sigma \in \text{aut}(\lie g)$, $m \in \mathbb{Z}_{>0}$ such that $\sigma^m= \text{Id}$ (i.e., $m$ is some multiple of $\left| \sigma \right|$; such an integer is called a \textit{period} of $\sigma$), and let $\zeta$ be a primitive $m^{th}$ root of unity.  Then $\lie g$ decomposes as a direct sum of $\sigma$--eigenspaces:
\[ \lie g = \bigoplus_{i=0}^{m-1} \lie g_{i}, \; \; \; \lie g_{i} = \left\{ x \in \lie g : \sigma(x) = \zeta^i x \right\},\]
We now define the \textit{loop algebra} $L(\lie g, \sigma, m)$ as 
\[ L(\lie g, \sigma, m) = \bigoplus_{s \in \mathbb{Z}} \lie g_{\overline{s}} \otimes t^s, \; \; \left[ x \otimes t^r, y \otimes t^s \right]= \left[ x, y\right] \otimes t^{r+s}\]
where $\overline{s}$ is the unique integer in $\left\{ 0, \ldots, m-1 \right\}$ modulo $m$.

In particular, when $\sigma = \mathrm{Id}$ and $m=1$, the resulting Lie algebra $L(\lie g, \mathrm{Id}, 1) = \lie g \otimes \mathbb{C}\left[ t^{\pm 1} \right]$ is often referred to as \textit{the} loop algebra (or sometimes the \textit{untwisted} loop algebra), and is usually denoted $L(\lie g)$.  If $\sigma = \mathrm{Id}$ and $m > 1$, then $L(\lie g, \mathrm{Id}, m) = \lie g \otimes \mathbb{C}\left[ t^{\pm m} \right]$, which is isomorphic as a Lie algebra to $L(\lie g, \mathrm{Id}, 1)$.  

\subsubsection{An alternate definition of the loop algebra $L(\lie g, \sigma, m)$.}\mbox{}

Let $\sigma \in \text{aut}(\lie g)$, $\sigma^m = \text{Id}$, and $\zeta$ a primitive $m^{th}$ root of unity.  Then we may extend $\sigma$ to a Lie algebra automorphism of the (untwisted) loop algebra by extending linearly the assignment
\begin{align*} 
\tilde{\sigma} : L(\lie g) & \rightarrow L(\lie g) \\  x \otimes t^k & \mapsto \zeta^{-k} \sigma(x) \otimes t^k, \end{align*}   
and define the loop algebra $L(\lie g, \sigma, m)$ as the set of fixed points of $L(\lie g)$ under the action of $\tilde{\sigma}$:
\[ L(\lie g, \sigma, m) = (L(\lie g))^\sigma = \left\{ a \in L(\lie g) \left. \right| \tilde{\sigma}(a) = a \right\}.\]
Then $L(\lie g, \sigma,m)$ is a Lie subalgebra of $L(\lie g)$, a relation which will be exploited in later sections.

\begin{lem}{\cite[Lemma 2.3] {ABP}}\label{period} Let $\sigma \in \text{aut}(\lie g)$, and suppose $\sigma^m = \Id$.  Then 
\[ L(\lie g, \sigma, m) \cong L(\lie g, \sigma, \left| \sigma \right|).\]
\end{lem}

\subsubsection{Ideals of loop algebras}\mbox{}

  Let $P \in \mathbb{C}\left[ t^m \right]$, and denote by $\left\langle P \right\rangle$ the ideal generated by $P$ in $\mathbb{C}\left[ t^{\pm m} \right]$.  Any ideal of a loop algebra $\twistg$ is given by such an ideal, in the following sense:
\begin{prop}[{\cite[\S 8.6, Lemma 8.6]{Kac}}]
\label{loopideal}
Let $I$ be an ideal of $\twistg$.  Then there exists a polynomial $P \in \mathbb{C}\left[ t^m \right]$ such that $I = \sum_{s \in \mathbb{Z}} \lie g_{\overline{s}} \otimes t^s\left\langle P \right\rangle$.  In particular any nonzero ideal of $\twistg$ is of finite codimension.
\end{prop}
\noindent Therefore any proper quotient of $\twistg$ is a finite--dimensional Lie algebra.  Following the notation of \cite{R}, for $P \in \mathbb{C}\left[ t^m \right]$ we set
\[ \twistg(P)  = \twistg  / \sum_{s \in \mathbb{Z}} \lie g_{\overline{s}} \otimes t^s\left\langle P \right\rangle. 
\] 

\subsection{Diagram automorphisms}\mbox{}

Let $\lie g = \lie g(X)$, $X = (X_{ij})_{1 \leq i, j \leq \ell}$ and $I_X = \left\{ 1, \ldots, \ell \right\}$.  A bijection $\tau: I_X \rightarrow I_X$ is a \textit{diagram automorphism} of $I_X$ if $X_{ij} = X_{\tau(i) \tau(j)}$.  Given such a diagram automorphism of $I_X$ and a Chevalley basis $\left\{ x^\pm_i, h_i\right\}_{1 \leq i \leq \ell}$ of $\lie g$, there exists a unique Lie algebra automorphism $\tau \in \text{aut}(\lie g)$ satisfying 
\begin{equation}\label{diagramauto} \tau(x^\pm_i) = x^\pm_{\tau(i)}.\end{equation}
We call $\sigma \in \text{aut}(\lie g)$ a \textit{diagram automorphism} of $\lie g$ if $\sigma$ is conjugate to some automorphism of the form (\ref{diagramauto}). 

There is a concrete description of these automorphisms for all types of finite--dimensional simple Lie algebras.  For example, when $\lie g$ is of type $A_n$ and $n> 1$ we have the following description.  There are two diagram automorphisms of $I_{A_n} = \left\{ 1, \ldots, n \right\}$ - the identity Id: $i \mapsto i$ and the `reflection' $\sigma: i \mapsto  n-i+1$ of order two.  Denote by $E_{i,j}$ the $(i,j)$ elementary matrix.  The Lie algebra automorphism $\sigma$ corresponding to this diagram automorphism has a very simple description when $\lie g$ is realized as the Lie algebra of traceless $n\times n$ matrices and the standard Chevalley basis 
\[ X^+_i = E_{i, i+1}, \; \; X^-_i = E_{i+1,i}, \; \; i =1, \ldots, n,\]
of $\lie g$ is chosen.  In this case $\sigma$ reflects all of the Chevalley generators $\left\{ E_{i, i+1}, E_{i+1, i} \right\}$ across the anti-diagonal:
\[ \sigma: \begin{cases} E_{i,i+1} &\longleftrightarrow \; \; E_{n-i+1, n-i+2}, \\ E_{i+1, i} &\longleftrightarrow \; \; E_{n-i+2 , n-i+1}. \end{cases}\]

The phrase `diagram automorphism' is used for bijections of the form (\ref{diagramauto}) because such a bijection provides a symmetry of the Dynkin diagram corresponding to $X$.  It is then easy to see that only for Lie algebras of type $A_n, D_n$ or $E_6$ do non--trivial diagram automorphisms exist, and that all of these are of order 2 with the exception of $D_4$, which has two diagram automorphisms of order 3 and 3 of order 2.  

\subsubsection{An example: $\lie g$ of type $A_2$; $\sigma$ the order 2 diagram automorphism of $A_2$}\mbox{}

  The order 2 diagram automorphism $\sigma$ of $A_2$ permutes the two nodes of the Dynkin diagram; i.e., $\sigma: 1 \leftrightarrow 2$.  Choose the standard Chevalley basis $X_1^+ = E_{1,2}, X_2^+ = E_{2,3}, X_1^- = E_{2,1}, X_2^- = E_{3,2}$ of $\lie g$.   Then under the action of $\sigma$, the vector space decomposition $\lie g = \lie g_0 \oplus \lie g_1$ into the two $\sigma$--eigenspaces $\lie g_0$ and $\lie g_1$  of $\lie g$ is  
\begin{align*}
&\lie g_0 = \text{span}_{\mathbb{C}}\left\{ \mat{0 & 0 & 0 \\ 1 & 0 & 0 \\ 0 & 1 & 0},  \mat{1 & 0 & 0 \\ 0 & 0 & 0 \\ 0 & 0 & -1} ,  \mat{0 & 1 & 0 \\ 0 & 0 & 1 \\ 0 & 0 & 0}  \right\}\\
& \lie g_1 = \text{span}_{\mathbb{C}} \left\{  \mat{0 & 0 & 0 \\ 0 & 0 & 0 \\ 1 & 0 & 0}, \mat{0 & 0 & 0 \\ 1 & 0 & 0 \\ 0 & -1 & 0}, \mat{1 & 0 & 0 \\ 0 & -2 & 0 \\ 0 & 0 & 1},  \mat{0 & 1  & 0 \\ 0 & 0 & -1 \\ 0 & 0 & 0},  \mat{0 & 0 & 1 \\ 0 & 0 & 0 \\ 0 & 0 & 0}\right\}.
\end{align*}

The twisted loop algebra $L(\lie g, \sigma, 2)$ is then
\[ L(\lie g, \sigma, 2) = \lie g_0 \otimes \mathbb{C}\left[ t^{\pm 2} \right] \oplus \lie g_1 \otimes t \mathbb{C}\left[ t^{\pm 2} \right]. \]

Although there are many more finite--order automorphisms of a Lie algebra $\lie g$ than the diagram automorphisms, these are convenient to work with because -- at least on the Chevalley generators of $\lie g$ -- their description is very simple.  Consequently, the twisted loop algebras $L(\lie g, \tau,m)$ corresponding to a diagram automorphism $\tau$ are easier to describe. Furthermore, we have the following proposition, due to Kac:

\begin{prop}[{\cite[\S 8.1, \S 8.3]{Kac}}]\label{kacprop}\mbox{}
\begin{enumerate}
\item[(i)] Let $\lie g$ be a simple finite--dimensional Lie algebra and $\sigma \in \aut(\lie g)$; $\sigma^m = \mathrm{Id}$.  Then there exists a Cartan subalgebra $\lie h$ of $\lie g$ such that $\sigma(\lie h)  = \lie h$, a triangular decomposition $\lie g = \lie n^- \oplus \lie h \oplus \lie n^+$, and an element $h \in (\lie h)^\sigma$ such that 
\[ \sigma = \mu \exp\left(\ad\left(\frac{2\pi i}{m}h\right)\right),\]
where $\mu$ is a diagram automorphism of $\lie g$ and $\mu(\lie n^\pm) = \lie n^\pm$. 
\item[(ii)] If $\sigma$ is a diagram automorphism, then $\lie g_0$ is a simple Lie algebra.  
\end{enumerate}
\end{prop}
The Cartan subalgebra $\lie h$ given in this proposition is briefly described as follows:

\begin{enumerate}
\item Let $\lie g_0$ be the fixed--point subalgebra of $\lie g$ under the action of $\sigma$.  Choose a maximal abelian subalgebra $\lie h_0$ of $\lie g_0$.  
\item The $\lie g$--centralizer of $\lie h_0$, $Z_{\lie g}(\lie h_0)$, is a Cartan subalgebra of $\lie g$.  Furthermore, this Cartan subalgebra contains some regular element $h \in (\lie h_0)^\sigma$.
\item The $\lie g$--centralizer of \textit{this} regular element is also a Cartan subalgebra.  We set $\lie h = Z_{\lie g}(h)$.  
\item With respect to the root space decomposition of $\lie g$ under the adjoint action of $\lie h$, there exists a choice of positive roots that are permuted by $\sigma$.  We set $\mu$ equal to the diagram automorphism corresponding to this permutation. 
\end{enumerate}

While the construction in Section \ref{loopdef} provides us with a loop algebra for any finite--order automorphism of $\lie g$, a theorem of Kac guarantees that, up to isomorphism, the loop algebras $L(\lie g, \mu, m)$ for a diagram automorphism $\mu$ exhaust all possibilities: 

\begin{thm}[{\cite[Proposition 8.5]{Kac}}]\label{Kac}
Let $\lie g$ be a simple finite--dimensional Lie algebra and $\sigma \in \aut(\lie g)$, $\sigma^m = \mathrm{Id}$. Then there exists a Lie algebra isomorphism $\Phi : L(\lie g, \mu, m) \rightarrow L(\lie g, \sigma, m)$ for some diagram automorphism $\mu \in \text{aut}(\lie g)$.
\end{thm}

The isomorphism $\Phi$ is defined as follows.  Let $h$ be the regular element given in the step (2) following Proposition \ref{kacprop}.  If $\alpha$ is an $\lie h$--root of $\lie g$, it follows that $\alpha(h) \in \mathbb{Z}$, and for $x \in \lie g_{\alpha}$, we set
\[ \Phi: x \otimes t^s \mapsto x \otimes t^{s + \alpha(h)}.\]
We illustrate all of the above with a simple example.  Let $\lie g = \lie{sl}_2$, and $\omega \in \text{aut}(\lie g)$ the involution of $\lie g$ defined by 
\[ \omega: X \mapsto -X^T\]
(with respect to the standard Chevalley basis $\left\{ \mat{ 0 & 0 \\ 1 & 0}, \mat{1 & 0 \\ 0 & -1}, \mat{0 & 1 \\ 0 & 0} \right\}$, $\omega$ is the \textit{Chevalley involution} of $\lie g$).  Then $L(\lie g, \omega, 2) \cong L(\lie g, \mu, 2)$ for some diagram automorphism $\mu$ of $\lie g$ with period 2.  Since the Dynkin diagram of $\lie g$ has only a single node, it has only the identity diagram automorphism $\mu = \mathrm{Id}$, and therefore 
\[ L(\lie g, \mu, 2) = \lie g \otimes \mathbb{C}\left[ t ^{\pm 2} \right] \cong \lie g \otimes \mathbb{C}\left[ t ^{\pm 1} \right].\]

The fixed--point subalgebra $\lie g_0$ of $\lie g$ under the action of $\omega$ is 
\[ \lie g_0 = \displaystyle{\left\{\left[ \begin{array}{cc}0&a \\ -a&0 \end{array} \right] :  a \in \mathbb{C} \right\} },\]
and a vector space basis of $\lie g_1$ is $\left\{ \mat{i & 1 \\ 1 & -i}, \mat{-i & 1 \\ 1 & i} \right\}$.  Therefore 
\[ L(\lie g, \omega, 2) = \text{span}_{\mathbb{C}} \left\{ \left[ \begin{array}{cc}0&1 \\ -1&0 \end{array} \right] \right\} \otimes \mathbb{C}\left[ t^{\pm 2}\right] \oplus \text{span}_{\mathbb{C}} \left\{ \mat{i & 1 \\ 1 & -i}, \mat{-i & 1 \\ 1 & i} \right\} \otimes t \mathbb{C}\left[ t^{\pm 2} \right]. \]

Since $\lie g_0$ is abelian, we have $\lie h_0 = \lie g_0$, and the corresponding Cartan subalgebra $Z_{\lie g}(\lie h_0) = \lie g_0$ as well (as can be verified by direct computation).  A triangular decomposition of $\lie g$ corresponding to this Cartan subalgebra is 
\[ \lie g = \mathbb{C}\left[ \begin{array}{cc}i & 1  \\ 1 &-i \end{array} \right] \oplus \mathbb{C}\left[ \begin{array}{cc} 0 & -i  \\ i & 0  \end{array} \right] \oplus \mathbb{C} \left[ \begin{array}{cc} -i & 1  \\ 1 & i \end{array} \right]\] 
(and these matrices provide us with a Chevalley basis of $\lie g$). 
The regular element $h\in Z_{\lie g}(\lie h_0)$ given in the proposition, and used in the construction of $\Phi$, must lie in $\mathbb{C}\mat{ 0 & 1 \\ -1 & 0}$, and must satisfy 
\[ -X^T = \omega(X) = \exp\left(\text{ad}\left(\frac{2\pi i}{m}h\right)\right)(X) \]
for any $X \in \lie g$.  Using this identity we find $h = \mat{0 & i / 2 \\ -i / 2 & 0}$.  

We now use this element $h$ to construct the isomorphism $\Phi : L(\lie g, \mathrm{Id}, 2) \rightarrow L(\lie g, \omega, 2)$, as follows.  If $\eta$ is a root of $\lie g$ and $x \in \lie g_{\eta}$, we define\footnote{If $\eta$ is a root of $\lie g$ (with respect to $Z_{\lie g}(h)$), it follows that $\eta(h) \in \mathbb{Z}$, so that this map is always well--defined.}
\[ \Phi: x \otimes t^s \mapsto x \otimes t^{s + \eta(h)}.\]
and extend $\Phi$ by linearity.

Let us denote the positive root of $\lie g$ by $\alpha$; then we have 
\[ \alpha\left( \mat{0 & i \\ -i & 0}\right)= 2, \text{ and therefore } \alpha \left(  \mat{0 & i/2 \\ -i/2 & 0}\right) = 1.\]
Now we can describe the action of $\Phi$ on a vector space basis of $L(\lie g, \mathrm{Id}, 2)$: for any $s \in \mathbb{Z}$, 
\[ \Phi: \begin{cases} \mat{i & 1 \\ 1 & -i } \otimes t^{2s} &\mapsto  \mat{i & 1 \\ 1 & -i } \otimes t^{2s+ (-\alpha)(h)} =  \mat{i & 1 \\ 1 & -i } \otimes t^{2s-1}, \\
\\
\mat{0 & 1 \\ -1 & 0}\otimes t^{2s} & \mapsto  \mat{0 & 1 \\ -1 & 0}\otimes t^{2s + (0)(h)} =   \mat{0 & 1 \\ -1 & 0}\otimes t^{2s} , \\
\\
\mat{-i & 1 \\ 1 & i } \otimes t^{2s} &\mapsto  \mat{-i & 1 \\ 1 & i } \otimes t^{2s+\alpha(h)}=  \mat{-i & 1 \\ 1 & i } \otimes t^{2s+1}. \end{cases}\]

The point of all of this is that $L(\lie g, \mathrm{Id}, 2) \cong \lie g \otimes \mathbb{C}\left[ t^{\pm 1} \right]$ is a simpler object to study than $L(\lie g, \omega, 2)$.  In particular, when considering the representation theory of $L(\lie{sl}_2, \sigma, m)$ for any finite--order automorphism $\sigma$ of $\lie{sl}_2$, it is sufficient to study only those representations of the `untwisted' $\lie{sl}_2 \otimes \mathbb{C}\left[ t^{\pm 1} \right]$.  

More generally, we may assume that, as far as the representation theory of $L(\lie g, \sigma, m)$ is concerned, we may assume that $\sigma$ is a diagram automorphism, and furthermore by Lemma \ref{period} that $m = \left| \sigma \right|$.  We will make these assumptions without further mention in the next two sections, and use the abbreviated notation $L(\lie g)^\sigma = L(\lie g, \sigma, m)$.  
 
\section{Representation theory}\label{representations}

In the following, all representations of a Lie algebra $\lie A$ are left $\lie A$--modules.  If $\lie{A}$ is a finite--dimensional Lie algebra, $\lie h$ a Cartan subalgebra of $\lie A$, $V$ a $\lie A$--module and $\lambda \in \lie h^*$, we define the \textit{$\lambda$--weight space $V_\lambda$ of $V$} as 
\[ V_{\lambda} = \left\{ v \in V : h.v = \lambda(h)v \right\}.\]

\subsection{Category--theoretic preliminaries}\mbox{}

Let $\mathcal M$ be the category whose objects are finite--dimensional $\lie A$--modules, and whose morphisms are $\lie{A}$--module homomorphisms. 
In this setting all objects have finite length, and we have the well-known 

\begin{theorem}[\textbf{Jordan--H\"{o}lder}]\mbox{}  Let $X$ be an object in $\mathcal M$.  Then 
\begin{enumerate}
\item [(i)] Any filtration of $X$ 
\begin{equation}\label{filtration}  \left\{ 0 \right\} = X_0 \subset X_1 \subset \cdots \subset X_\ell = X\end{equation} such that all $X_i / X_{i-1}$ are nonzero can be refined to a Jordan-H\"{o}lder series, i.e. a filtration 
\begin{equation}\label{JH}\left\{ 0 \right\} = X'_0 \subset X'_1 \subset \cdots \subset X'_m = X\end{equation}
such that all quotients $X'_i / X'_{i-1}$ are simple.
\item [(ii)] Any two Jordan-H\"{o}lder series of $X$ are unique (i.e., determine the same set of simple objects $X'_i /X'_{i-1}$ up to isomorphism).
\end{enumerate}
\end{theorem}

Given any filtration of an object $X \in \mathcal M$ as in (\ref{filtration}), we will refer to a quotient $X_i / X_{i-1}$ of a submodule $X_i$ pf $X$ as a \textit{subquotient} of $X$, and to the irreducible subquotients occurring in the Jordan-H\"{o}lder series (\ref{JH}) as the \textit{constituents} of $X$.

\begin{defn} 
\mbox{}
\begin{enumerate}
\item [(i)] Two indecomposable objects $V_1, V_2 \in \mathcal M$ are \textit{linked}, written $V_1 \sim V_2$, if there do not exist Abelian subcategories $\mathcal M_1, \mathcal M_2$ such that $\mathcal M \cong \mathcal M_1 \oplus \mathcal M_2$ with $V_1 \in \mathcal M_1$, $V_2 \in \mathcal M_2$.  More generally, two objects $U, V \in \mathcal M$ are linked if every indecomposable summand of $U$ is linked to every indecomposable summand of $V$.  We will say that an object $V \in \mathcal{M}$ is linked if there exists some object $W$ in $\mathcal{M}$ such that $V \sim W$.  The relation $\sim$, when restricted to the collection of linked objects\footnote{without this restriction, the relation $\sim$ is reflexive and symmetric, but not transitive.  For example, if $W_1$ and $W_2$ are two objects in $\mathcal{M}$ which are \textit{not} linked, then $W = W_1 \oplus W_2$ is not linked to itself or anything else in $\mathcal{M}$.}, is an equivalence relation.   
\item [(ii)] A block of $\mathcal M$ is an equivalence class of linked objects.
\end{enumerate}
\end{defn}

An abelian category $\mathcal C$ is \textit{indecomposable} if, for any equivalence of categories $\mathcal C \cong \mathcal C_1 \oplus \mathcal C_2$, we have either $\mathcal C \cong \mathcal C_1$ or $\mathcal C \cong \mathcal C_2$. 

\begin{prop}[{\cite[Proposition 1.1] {EM}} \label{catdecomp}]
The category $\mathcal M$ admits a unique decomposition into a direct sum of indecomposable abelian subcategories: $\mathcal M = \sum_{\alpha \in I}\mathcal M_{\alpha}$. 
\end{prop}
The indecomposable subcategories which occur in this decomposition are the blocks of the category.    
\subsection{Representations of finite--dimensional simple complex Lie algebras}\label{greps}\mbox{}

We briefly recall here some of the standard theory concerning the representations of a simple finite--dimensional Lie algebra $\lie g$.  While these definitions and results are well known (see \cite{H}, for example), we repeat them here to motivate similar constructions for the loop algebras.  

\subsubsection{Highest--weight modules of $\lie g$}\mbox{}

  A $\lie g$--module $V$ is called a \textit{highest--weight module with highest weight $\lambda \in \lie h^*$} (or a $\lambda$--highest weight module, or sometimes just a highest--weight module) if there exists some $v \in V$ with $U(\lie g).v = V$ such that 
               \[\lie n^+.v = 0, \; \; \; h.v = \lambda(h)v \; \text{ for all } h \in \lie h.\]

\subsubsection{Verma modules of $\lie g$}\mbox{}

 Let $\lambda \in \lie h^*$, and define an action of the subalgebra $\lie h \oplus \lie n^+ \subseteq \lie g$ on a one--dimensional complex vector space $\mathbb{C}_{\lambda}=\mathbb{C}1$ as follows:
\begin{align*}
&\lie n^+ . 1 = 0, \\
&h. 1 = \lambda(h)1, \; \text{ for all } h \in \lie h.
\end{align*}
Now define the \textit{Verma module $Z(\lambda)$ with highest weight $\lambda$} as 
\[ Z(\lambda) := U(\lie g) \otimes_{U(\lie h \oplus \lie n^+)}\mathbb{C}_\lambda.\]
Then $Z(\lambda)$ has a natural left $\lie g$--module structure.  $Z(\lambda)$ may also be defined as the quotient space $U(\lie g) / I_{\lambda}$, where $I_\lambda$ is the left ideal of $U(\lie g)$ generated by $\lie n^+, (h - \lambda(h))1$ for all $h \in \lie h$.   $\lie g$ then acts by left multiplication on the left cosets of $Z(\lambda)$.  

\begin{thm}\mbox{}
\begin{enumerate}
\item[(i)] Verma modules are `universal highest--weight modules': $Z(\lambda)$ is a $\lambda$--highest weight module, and if $W$ is any $\lambda$--highest weight module generated by $w \in W_{\lambda}$ then the map $1 \mapsto w$ (where $1$ is the left coset of the identity $1 \in U(\lie g)$) extends to a surjective homomorphism $Z(\lambda) \twoheadrightarrow W$.  
\item[(ii)] The collection of all Verma modules is indexed (up to isomorphism) by the space of linear forms $\lie h^*$, via the bijection $\lambda \leftrightarrow Z(\lambda)$.  
\item[(iii)]  Any Verma module $Z(\lambda)$ has a unique maximal submodule and hence a unique irreducible quotient, which we denote by $V(\lambda)$.  
\end{enumerate}
\end{thm}

\subsubsection{Irreducible finite--dimensional $\lie g$--modules}\mbox{}

As noted above, any Verma module $Z(\lambda)$ has a unique irreducible quotient $V(\lambda)$.  But this irreducible quotient is not necessarily finite--dimensional.  The necessary and sufficient condition for this is given by the following  

\begin{thm}\mbox{}
\begin{enumerate}
\item Let $Z(\lambda)$ be the Verma module with highest weight $\lambda \in \lie h^*$ and $V(\lambda)$ its irreducible quotient.  Then $V(\lambda)$ is finite--dimensional if and only if $\lambda \in P^+$.  
\item Let $V$ be a finite--dimensional irreducible $\lie g$--module.  Then $V$ is $\lambda$--highest weight for some $\lambda \in P^+$ (hence $V$ is isomorphic to the irreducible quotient of $Z(\lambda)$). 
\item \textbf{Weyl's Theorem}: Let $W$ be a finite--dimensional $\lie g$--module.  Then $W$ can be written (uniquely, up to permutation of summands) as a direct sum 
\[ W \cong \bigoplus_{i=1}^\ell V(\lambda_i)\]
of irreducible $\lie g$--modules $V(\lambda_i)$ .
\end{enumerate}
\end{thm}

\subsubsection{The category \mbox{$\mathcal{C}$}}\mbox{}

  Let $\mathcal{C}$ be the category of finite--dimensional representations of $\lie g$.  This category $\mathcal{C}$ is \textit{semisimple} -- any object in $\mathcal{C}$ is completely reducible (this is Weyl's theorem).  This, in turn, is equivalent to the statement that the blocks of the category $\mathcal{C}$ are in bijective correspondence with the isomorphism classes of finite--dimensional irreducible modules in $\mathcal{C}$:
\begin{equation}\label{gblocks}
 \mathcal C = \bigoplus_{\eta \in P^+} \mathcal C_\eta\end{equation}
where, for $\eta \in P^+$, $\mathcal{C}_\eta$ is the subcategory of $\mathcal{C}$ which consists of all direct sums $V(\eta)^{\oplus i}$, $i \in \mathbb{N}$.  

\subsection{Representations of loop algebras}\label{reps}\mbox{}

We would now like to exploit the concepts of a highest--weight representation and Verma module in a different setting; namely, that of the finite--dimensional representation theory of a loop algebra $L(\lie g)^\sigma$.  When we move from $\lie g$ to $L(\lie g)^\sigma$, we evidently lose two properties which $\lie g$ enjoyed: that of being finite--dimensional and that of being semisimple.  If we restrict our attention to those representations of the loop algebra which are finite--dimensional, both of these properties can be regained, in a certain sense.  

The definition of a highest--weight representation of $\lie g$ above depended entirely upon the triangular decomposition of $\lie g = \lie n^- \oplus \lie h \oplus \lie n^+$.  We will use a similar decomposition for a loop algebra.  In the case of the untwisted loop algebra $L(\lie g) = \lie g \otimes \mathbb{C}\left[ t^{\pm 1} \right]$, we have the decomposition 
\[ L(\lie g) = L(\lie n^-) \oplus L(\lie h) \oplus L(\lie n^+), \]
where $L(\lie n^\pm), L(\lie h)$ are the loop algebras of the subalgebras $\lie n^\pm, \lie h$.  If the triangular summands $\lie n^\pm, \lie h$ of $\lie g$ remain invariant under the action of $\sigma$, then $\sigma$ is an automorphism of these Lie subalgebras and we can define the twisted loop (sub)algebras $L(\lie n^\pm, \sigma, m), L(\lie h, \sigma, m)$, and we have 
\begin{equation}\label{twistedtriangle} L(\lie g)^\sigma = L(\lie n^-, \sigma, m) \oplus L(\lie h, \sigma, m) \oplus L(\lie n^+, \sigma, m).\end{equation}
Since $\sigma$ is a diagram automorphism of $\lie g$, we have $\sigma(\lie n^\pm) \subseteq \lie n^\pm, \; \sigma(\lie h) \subseteq \lie h$, and so we have the decomposition (5), in which case we will write  $L(\lie n^\pm)^\sigma = L(\lie n^\pm, \sigma, m)$, and $L(\lie h)^\sigma = L(\lie h, \sigma, m)$. 


\subsubsection{Loop--highest weight representations} \mbox{}

\begin{defn}
Let $V$ be an $L(\lie g)^\sigma$--module and $\Lambda \in (L(\lie h)^\sigma)^*$.  We say $V$ is a \textit{loop highest--weight module with highest weight $\Lambda$} if there exists an element $v \in V$ with $U(L(\lie g)^\sigma).v = V$, such that 
\[ L(\lie n^+)^\sigma.v = 0, \; \; \; h.v = \Lambda(h)v \; \text{ for all } h \in L(\lie h)^\sigma.\]
\end{defn}
Of course, this definition closely parallels that of the highest--weight representations of a simple finite--dimensional Lie algebra $\lie g$. 

\begin{defn}
Let $\Lambda \in (L(\lie h)^\sigma)^*$, and let $\mathbb{C}_\Lambda$ be the 1--dimensional representation of the subalgebra $L(\lie h)^\sigma \oplus L(\lie n^+)^\sigma \subseteq L(\lie g)^\sigma$ defined by 
\begin{align*}
x.1 &= 0,  \; \; \; \;  x \in L(\lie n^+)^\sigma \\
h. 1 &= \Lambda(h),  \; \; \; \; h \in L(\lie h)^\sigma.
\end{align*}
Then we set
\[ Z(\Lambda) = U(L(\lie g)^\sigma)\otimes_{U(L(\lie h)^\sigma \oplus L(\lie n^+)^\sigma)}\mathbb{C}_\Lambda.\]
\end{defn}
We could call $Z(\Lambda)$ the \textit{loop Verma module} with weight $\Lambda$, although this terminology has not been used in the literature.

\subsubsection{Restrictions and evaluation representations}\mbox{}

Recall that $\lie g_0 = \lie g^\sigma$, the fixed--point subalgebra of $\lie g$ under the automorphsim $\sigma: \lie g \rightarrow \lie g$.  Because of the algebra inclusions 
\[\xymatrix@R=0.03in{ \lie g_0 \ar@{^{(}->}[r] &  L(\lie g)^\sigma \ar@{^{(}->}[r] & L(\lie g), \\
                   x \ar@{|->}[r] & x \otimes 1 \ar@{|->}[r]&  x \otimes 1 ,}\]
 any representation $W$ of $L(\lie g)$ can be restricted to yield a representation of $L(\lie g)^\sigma$ or of the simple Lie algebra $\lie g_0$ (in particular, of $\lie g$ when $\sigma = \mathrm{Id}$). 

But we also have a method of extending a representation of $\lie g$ to a representation of $L(\lie g)^\sigma$.  Given any $\ell$--tuple $\vec{a} = (a_1, \ldots, a_\ell) \in \mathbb{C}^\ell$ there exists a Lie algebra homomorphism 
\begin{align*}
\ev_\vec{a} : L(\lie g)^\sigma & \rightarrow \bigoplus_{i = 1}^\ell \lie g\\
x \otimes t^s &\mapsto (a_1^s x, \ldots, a_\ell^s x).
\end{align*}
We call $\ev_{\vec{a}}$ the \textit{evaluation homomorphism} (at $\vec{a}$).  Let $V_1, \ldots, V_\ell$ be $\lie g$--modules, via the homomorphisms $\lie g \stackrel{\psi_i}{\longrightarrow} \text{End}_{\mathbb{C}}(V_i)$.  
Then their tensor product $\otimes_{i=1}^\ell V_i$ is a representation of $\bigoplus_{i=1}^\ell \lie g$ via the Lie algebra homomorphism
\begin{align*} \otimes \psi_i: \oplus \lie g &\longrightarrow  \text{End}_{\mathbb{C}}\left(\bigotimes_{i=1}^\ell V_i\right) \\
                   (x_1, \ldots, x_\ell) &\mapsto \sum_{i=1}^\ell \mathrm{Id} \otimes \cdots \otimes \mathrm{Id} \otimes \psi_i(x_i) \otimes \mathrm{Id} \otimes \cdots \otimes \mathrm{Id}.\end{align*}
By composing this homomorphism with an evaluation homomorphism we obtain a representation of $L(\lie g)^\sigma$:
\[ L(\lie g)^\sigma \stackrel{\ev_{\vec{a}}}{\longrightarrow} \bigoplus_{i=1}^\ell \lie g \stackrel{\otimes \psi_i}{\longrightarrow } \text{End}_{\mathbb{C}}\left(\bigotimes_{i=1}^\ell V_i\right).\]

Henceforth we shall be concerned only with the finite--dimensional representations of $L(\lie g)^\sigma$.  Thus we define $\mathcal{F}^\sigma$ to be the category whose objects are finite--dimensional $L(\lie g)^\sigma$--modules, and whose morphisms are $L(\lie g)^\sigma$--module homomorphisms.  The evaluation representations (of finite--dimensional $\lie g$--modules) described above already give us many examples of objects in $\mathcal{F}^\sigma$.  

\subsubsection{The Weyl modules}\mbox{} 

In the category $\mathcal{F}^\sigma$ there exist certain maximal loop--highest weight modules; these are the \textit{Weyl modules} of $L(\lie g)^\sigma$.  

When $\sigma \in \aut(\lie g)$ is a diagram automorphism, the subalgebra $\lie g_0$ is a simple finite--dimensional Lie algebra (see {\cite[\S 8]{Kac}}).  Recall that we have set $n = \text{rank}(\lie g)$, and $I = \left\{ 1, \ldots, n \right\}$.

Let $\text{rank}(\lie g_0)= n_0$, let $I_0 = \left\{ 1, \ldots, n_0 \right\}$ and denote the simple roots of ($\lie g_0, \lie h_0$) by $\tilde{\alpha}_1, \ldots, \tilde{\alpha}_{n_0}$.  There exists a natural identification $I_0 \rightarrow I / \sigma$ of $I_0$ with the  $\sigma$--orbits of $I$ (see \cite{Lus} for details), and 
we choose an identification
\begin{align*}
\iota: I_0 & \hookrightarrow I, \\ i &\mapsto \iota(i),
\end{align*}
of each $i \in I_0$ with a representative of an orbit in $I$ which lifts this identification; i.e., an injection $\iota$ satisfying 
\[ \xymatrix@R=0.03in{ & I \ar@{->>}[dd] \\ I_0 \ar@{^(->}^{\iota}[ur] \ar[dr] & \\ & I / \sigma}\]
 For each $i \in I_0$, we then define $S(i) := \left| \left\langle \sigma \right\rangle_{\iota(i)} \right|$, where $\left\langle \sigma \right\rangle_{\iota(i)} $ is the stabilizer of $\iota(i) \in I$ under the action of $\sigma$. \\
For example, when $\lie g$ is of type $E_6$, there is one order 2 diagram automorphism $\sigma \in \aut(\lie g)$.  The fixed--point subalgebra $\lie g_0$ of $\lie g$ is of type $F_4$.  One choice of embedding $\iota: I_0 \hookrightarrow I$ is pictured here (the only other choice sends $2 \mapsto 4$ and $1\mapsto 5$).

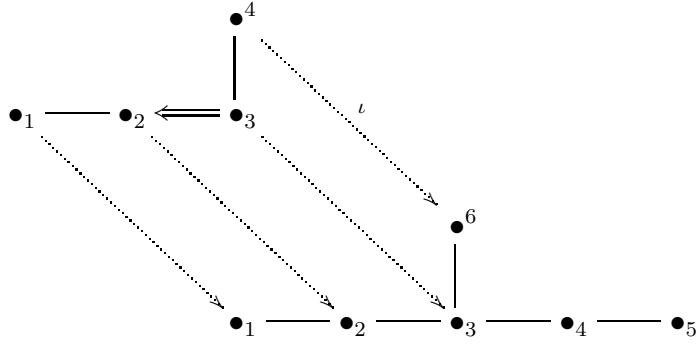
\[ \xymatrix{ & & \bullet^4 \ar@<-0.4ex>@{.>}[ddrr]^\iota \\ \bullet_1 \ar@{.>}[ddrr] \ar@<0.4ex>@{-}[r] &  \bullet_2 \ar@{.>}[ddrr] \ar@<0.4ex>@{<=}[r] &  \bullet_3 \ar@{.>}[ddrr] \ar@<0.6ex>@{-}[u] \\
& & & & \bullet^6\\
& & \bullet_1  \ar@<0.4ex>@{-}[r] &  \bullet_2 \ar@<0.4ex>@{-}[r] &  \bullet_3  \ar@<0.6ex>@{-}[u] \ar@<0.4ex>@{-}[r] & \bullet_4 \ar@<0.4ex>@{-}[r] & \bullet_5} \]
In this case we have $S(1) = S(2) = 1$ and $S(3) = S(4) = 2$.

Recall that $m$ is the order of $\sigma$, and $\lie g_\epsilon$ is the $\zeta^\epsilon$--eigenspace of $\lie g$ under the action of $\sigma$.  For $\epsilon \in \left\{ 0, \ldots, m-1 \right\}$, $i \in I_0$ and $y_j \in \left\{x_j^\pm, h_j \right\}$, we define the elements 
\begin{align*}
y_{i, \epsilon} = \sum_{j=0}^{m-1}(\zeta^\epsilon)^j\sigma^j(y_{\iota(i)}).
\end{align*}
The subspace $\lie h_\epsilon = \lie g_\epsilon \cap \lie h$ is then spanned by $\left\{ h_{i, \epsilon} \right\}_{ i \in I_0}$. 

\begin{defn}
Let $\mathcal{P}$ be the set of all $n_0$--tuples of polynomials $\bpi^\sigma = \left( \pi_1(u), \ldots, \pi_{n_0}(u) \right) \in (\mathbb{C}\left[ u \right])^{n_0}$ such that $\pi_i(0) = 1$ for $i = 1, \ldots, n_0$.  

\end{defn}
We will write an element $\bpi^\sigma \in \mathcal{P}$ as $\bpi^\sigma = \left( \pi_i \right)_{i \in I_0}$.  The set $\mathcal{P}$ is a monoid under pointwise multiplication of polynomials: 
\[ \left( \pi_1, \ldots, \pi_{n_0} \right)\left( \pi_1', \ldots, \pi_{n_0}' \right) = \left( \pi_1 \pi_1', \ldots, \pi_{n_0}\pi_{n_0}' \right).\]
For $i \in I_0$ and $a \in \mathbb{C}^\times$ we define  
\[ \bpi^\sigma_{i,a} = \left( (1-a^{S(i)} u)^{\delta_{ij}} \right) \in \mathcal{P}. \]
Then any $\bpi^\sigma$ can be written as a product 
\[ \bpi^\sigma  = \prod_{j=1}^\ell \prod_{i = 1}^{n_0}(\bpi^\sigma_{i,a_j})^{s_{ij}}, \; \; \; a_j^m \neq a_k^m \text{ for } j \neq k,\]
where $s_{ij}$ is the multiplicity of the root $a_j^{-S(i)}$ in the $i^{th}$ component $\pi_i(u)$ of $\bpi^\sigma$.  

If $\sigma = \text{Id}$, we will write $\bpi^\sigma = \bpi$. We define a function on the elements $\bpi_{i, a} \in \mathcal{P}^{\Id}$ as follows:
\begin{equation}\label{boldr} \mathbf{r} : \bpi_{i, a} \mapsto \bpi^\sigma_{\sigma^\epsilon(i), \zeta^\epsilon a}\end{equation}
where $\epsilon$ is the unique element in $ \left\{ 0, \ldots, m-1 \right\}$ such that $\sigma^{\epsilon}(i) \in I_0$.  Since $\mathcal{P}^{\Id}$ is multiplicatively generated by $\left\{ \bpi_{i, a} \right\}_{i \in I, a \in \mathbb{C}^\times}$, there exists a unique monoid homomorphism $\mathbf{r}: \mathcal{P}^{\Id} \twoheadrightarrow \mathcal{P}$ satisfying (\ref{boldr}). 

To any $\bpi \in \mathcal{P}^{\Id}$ we can associate a list of dominant integral weights (with respect to $(\lie g, \lie h)$), described as follows: let $\bpi = \prod_{j=1}^\ell \prod_{i\in I} (\bpi_{i,a_j})^{s_{ij}}$, and define 
\[ \lambda_j = \sum_{i \in I} s_{ij} \omega_i, \; \; j = 1, \ldots, \ell.\]
If, for $\lambda \in P^+$ and $a \in \mathbb{C}^\times$ we define 
\[ \bpi_{\lambda, a} = \left( (1 -au)^{\lambda(h_i)} \right)_{i \in I},\]
then any $\bpi \in \mathcal{P}^{\Id}$ has a unique factorization 
\[ \bpi = \prod_{j=1}^\ell \bpi_{\lambda_j, a_j}, \; \; \; a_j \neq a_k \text{ for } j \neq k.\]

\begin{defn}
Let $\bpi^\sigma = \prod_{j=1}^\ell \prod_{i = 1}^{n_0}(\bpi^\sigma_{i,a_j})^{s_{ij}} \in \mathcal{P}$, and let $W(\bpi^\sigma)$ be the quotient  
\[ W(\bpi^\sigma) = U(L(\lie g)^\sigma) / J_{\bpi^\sigma},\] 
where $J_{\bpi^\sigma}$ is the left ideal in $U(L(\lie g)^\sigma)$ generated as follows: 
\[ J_{\bpi^\sigma} = \left\langle L(\lie n^+)^\sigma, \; \; \; \left(x_{i,0}^- \otimes 1\right)^{\deg(\pi_i)+1}, \; \; \; h_{i,\epsilon}\otimes t^{ms - \epsilon} - \sum_{j=1}^\ell a_j^{ms-\epsilon}s_{ij} \right\rangle \]
for $i \in I_0$, $\epsilon \in \left\{ 0, \ldots, m-1\right\}$, and $s \in \mathbb{Z}$.  $L(\lie g)^\sigma$ acts on $W(\bpi^\sigma)$ by left multiplication of the cosets of $W(\bpi^\sigma)$.  
\end{defn}
The inclusion of $L(\lie n^+)^\sigma$,  $ h_{i,\epsilon}\otimes t^{ms - \epsilon} - \sum_{j=1}^\ell a_j^{ms-\epsilon}s_{ij}$ in $J_{\bpi^\sigma}$ ensures that $W(\bpi^\sigma)$ is loop--highest weight, while $\left(x_{i,0}^- \otimes 1\right)^{\deg(\pi_i)+1} \in J_{\bpi^\sigma}$ ensures that $L(\lie n^-)^\sigma$ acts locally nilpotently on $W(\bpi^\sigma)$.  

\section{Classification of Weyl modules, simple objects and blocks}\label{loopclassification}

\subsection{Parametrization of Weyl modules}\mbox{}

It is evident from the definition above that the modules $W(\bpi^\sigma)$ are loop--highest weight for all $\bpi^\sigma \in \mathcal{P}$.  But they are also maximal among all such finite--dimensional $L(\lie g)^\sigma$--modules.

\begin{thm}[{\cite[Theorem 1, Proposition 2.1]{CPweyl}, {\cite[Theorem 2, Proposition 4.2]{CFS}}}] \mbox{}
\begin{enumerate}
\item[(i)] $\dim(W(\bpi^\sigma)) < \infty$ for all $\bpi^\sigma \in \mathcal{P}$.  
\item[(ii)] Let $V$ be a loop--highest weight module in $\mathcal{F}^\sigma$.  Then there exists a unique $\bpi^\sigma \in \mathcal{P}$ such that $V$ is an $L(\lie g)^\sigma$--module quotient of $W(\bpi^\sigma)$.
\item[(iii)] $W(\bpi^\sigma)$ has a unique irreducible quotient.
\item[(iv)] Let $\bpi^\sigma \in \mathcal{P}$.  Then there exists $\bpi \in \mathbf{r}^{-1}(\bpi^\sigma)$ such that, as $L(\lie g)^\sigma$--modules, 
\[ W(\bpi^\sigma) \cong W(\bpi)\left.\right|_{L(\lie g)^\sigma},\]
and this statement is also true if we replace the Weyl module W(\bpi) with its irreducible quotient $V(\bpi)$.
\end{enumerate}
\end{thm}
Therefore the isomorphism classes of simple objects in $\mathcal{F}^\sigma$ are, along with the Weyl modules, parametrized by $\mathcal{P}$.  
We denote the unique irreducible quotient of a Weyl module $W(\bpi^\sigma)$ by $V(\bpi^\sigma)$, and summarize these facts with the diagram
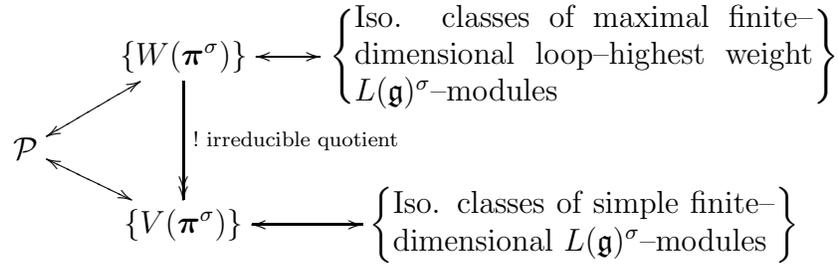
\[ \xymatrix@R=0.08in{ & \left\{ W(\bpi^\sigma) \right\} \ar@{<->}[r] \ar@{->>}[dd]^{!  \text{ irreducible quotient }} & \left\{\text{\parbox{2.4 in}{Iso. classes of maximal finite--dimensional loop--highest weight  $L(\lie g)^\sigma$--modules}}\right\}
\\
\mathcal{P} \ar@{<->}[ur] \ar@{<->}[dr] \\ & \left\{ V(\bpi^\sigma) \right\} \ar@{<->}[r]&  \left\{\text{\parbox{2 in}{Iso. classes of simple finite--dimensional $L(\lie g)^\sigma$--modules}}\right\}}\]
The `maximal' property here is that described by the universal property of the Weyl modules given in (ii) of the theorem. 

We can say more about the parametrization $\mathcal{P} \leftrightarrow \left\{ V(\bpi^\sigma) \right\}$: when $\sigma = \Id$ and $\bpi = \prod_{j=1}^\ell \bpi_{\lambda_i, a_i}$, $a_j \neq a_k$, the irreducible $L(\lie g)$--module $V(\bpi)$ is isomorphic to the evaluation module 
\begin{equation}\label{irreps} \xymatrix{ L(\lie g) \ar@{->>}[r]^-{\ev_{\vec{a}}}  & \bigoplus_{i=1}^\ell \lie g  \ar[r]^-{\otimes \psi_i} &  \text{End}_{\mathbb{C}}\left(\otimes_{i=1}^\ell V(\lambda_i)\right),}\end{equation}
where $\vec{a} = \left( a_1, \ldots, a_\ell\right)$, and $\psi_i : \lie g \rightarrow \text{End}_{\mathbb{C}}(V(\lambda_i))$ is the representation of $\lie g$ on the simple finite--dimensional $\lie g$--module $V(\lambda_i)$ with highest weight $\lambda_i$ (see \cite{C}, \cite{CPunitary}, \cite{B}, \cite{CPweyl}, and Sections \ref{multiloops} and \ref{additional} for more on this realization) . 

Although this provides us with a complete list of simple finite--dimensional $L(\lie g)^\sigma$--modules, the category $\mathcal{F}^\sigma$ is not semisimple -- the Weyl modules provide examples of indecomposable and (in general) reducible objects in $\mathcal{F}^\sigma$.  Hence the blocks of $\mathcal{F}^\sigma$ are not parametrized by $\mathcal{P}$, but by a collection of equivalence classes of $\mathcal{P}$ which we now describe.

\subsubsection{The spectral characters of $L(\lie g)^\sigma$}\mbox{}

\begin{defn}[{\cite{CM}}]

Let $\Xi = \left\{ \chi : \mathbb{C}^\times \rightarrow P/Q : \left| \supp(f) \right| < \infty \right\}$, and for $\lambda \in P$ let $\overline{\lambda} = \lambda + Q \in P/Q$.  Then $\Xi$ is an additive group, and any $\chi \in \Xi$ can be written uniquely as a sum
\[ \chi = \sum_{j=1}^\ell \chi_{\lambda_j, a_j}, \; \; \; a_j \neq a_k \text{ for }j \neq k,\; \; \text{ where } \chi_{\lambda, a}(z) = \begin{cases} \overline{\lambda}, & z = a \\ 0, & z \neq a.\end{cases}\]
\end{defn}
\noindent In \cite{CM}, the elements $\chi \in \Xi$ are called the \textit{spectral characters} of $L(\lie g)$. 

Define a map $\varpi$ on the elements $\bpi_{i, a} \in \mathcal{P}^{\Id}$ by setting 
\[ \varpi: \bpi_{i, a} \mapsto \chi_{\omega_i, a}.\]
Since $\mathcal{P}^{\Id}$ is (multiplicatively) generated by $\left\{ \bpi_{i, a} \right\}_{i \in I, a \in \mathbb{C}^\times}$ and $\Xi$ is (additively) generated by $\left\{ \chi_{\omega_i, a} \right\}_{i \in I, a \in \mathbb{C}^\times}$, $\varpi$ extends uniquely to a surjective monoid homomorphism $\varpi: \mathcal{P}^{\Id} \twoheadrightarrow \Xi$, and we set $\chi_{\bpi} = \varpi(\bpi)$.  

\begin{defn} Two spectral characters $\chi_1, \chi_2 \in \Xi$ are $\sigma$--\textit{equivalent}, written $\chi_1 \sim_\sigma \chi_2$, if there exist $\bpi_i \in \mathcal{P}^{\Id}$ such that $\chi_i = \chi_{\bpi_i}$ and $\mathbf{r}(\bpi_1) = \mathbf{r}(\bpi_2)$ .
\end{defn}

It is routine to show that the relation $\sim_\sigma$ is an equivalence relation on the group $\Xi$ (see {\cite[\S 3.3]{S}}).  If $\chi \in \Xi$, we denote the corresponding equivalence class by $\overline{\chi}$.  The binary operation $\overline{\chi_{\bpi_1}} + \overline{\chi_{\bpi_2}} = \overline{\chi_{\bpi_1\bpi_2}}$ is then well--defined.

\begin{defn}[{\cite{CM}, \cite{S}}] \mbox{}
\begin{enumerate}
\item [(i)]  The group of equivalence classes of $\sigma$--equivalent spectral characters is called the group of \textit{twisted spectral characters}, denoted $\Xi^\sigma$.  
\item [(ii)] Let $V \in \mathcal{F}^\sigma$ with constituents $\left\{ V(\bpi_i^\sigma) \right\}_{i=1, \ldots, s}$, and let $\bpi_i \in \mathbf{r}^{-1}(\bpi^\sigma_i) \subseteq \mathcal{P}^{\Id}$.  The module $V$ has \textit{twisted spectral character} $\overline{\chi}$ if $\overline{\chi_{\bpi_i}} = \overline{\chi}$ for $i = 1, \ldots, s$. 
\end{enumerate}
\end{defn}

Let $\mathcal{F}^\sigma_{\overline{\chi}}$ be the subcategory of $\mathcal{F}^\sigma$ consisting of all finite--dimensional $L(\lie g)^\sigma$--modules with twisted spectral character $\overline{\chi}$.  In the category $\mathcal{F}^\sigma$, we can show that
\begin{enumerate}
\item [(i)] Any indecomposable object in $\mathcal{F}^\sigma$ has a twisted spectral character, and 
\item[(ii)] Any two simple objects in $\mathcal{F}^\sigma$ with the same twisted spectral character
are linked.
\end{enumerate}
The proof of (i) uses the property that there are no non--trivial extensions between modules with distinct spectral characters (i.e., if $V_i \in \mathcal{F}^\sigma_{\overline{\chi_i}}$ for $i=1,2$ with $\overline{\chi_1} \neq \overline{\chi_2}$, then $\text{Ext}^1(V_1, V_2) = 0$; see {\cite[Lemma 5.2]{CM}} and {\cite[Lemma 3.25]{S}}).  If $\left\{ V_i \right\}$ is a list of constituents of an object $W$ and among these objects, there is one which has a spectral character which is distinct from the rest, we are able to use this property to produce a non--trivial decomposition of $W$. 

For (ii), we begin with two simple objects $V_1, V_2$ which share the same twisted spectral character and produce a nontrivial element $V \in \text{Ext}^1(V_1, V_2)$.  This nontrivial extension provides us with a chain of $L(\lie g)^\sigma$--homomorphisms between indecomposable modules
\[ \xymatrix{ V_2 \ar@{^(->}[r]  & V \ar@{->>}[r] & V_1,}\]
which implies that $V_1$ and $V_2$ are linked (the `linkage' provided by such a chain of maps between indecomposable objects is, in fact, equivalent to the category--theoretic linkage defined above; see  {\cite[Lemma 2.5]{CM}}).  Note that the construction of such a non--trivial extension between simple objects in the category $\mathcal{C}$ of finite--dimensional $\lie g$--modules is prohibited by Weyl's theorem -- in the semisimple category $\mathcal{C}$, all distinct non--isomorphic simple objects lie in distinct blocks. 
 
\noindent From properties (i) and (ii) follows the main theorem concerning block decomposition of $\mathcal{F}^\sigma$:
\begin{thm} [{\cite[Theorem 1]{CM}}, {\cite[Theorem 2]{S}}]
The blocks of $\mathcal{F}^\sigma$ are in bijective correspondence with the set $\Xi^\sigma$.  More precisely, we have a decomposition of $\mathcal{F}^\sigma$ into indecomposable subcategories 
\[ \mathcal{F}^\sigma = \bigoplus_{\overline{\chi} \in \Xi^\sigma} \mathcal{F}^\sigma_{\overline{\chi}}.\]
\end{thm}
\noindent Compare with the decomposition (\ref{gblocks}) given in Section \ref{greps}.  This theorem was proved in the untwisted $\sigma = \Id$ case in \cite{CM} and extended to the category $\mathcal{F}^\sigma$ for an arbitrary diagram automorphism $\sigma$ in \cite{S}.

\section{Representations of multiloops}\label{multiloops}
The description of simple objects in $\mathcal{F}^\sigma$ as tensor products of evaluation representations (as described above, see (\ref{irreps})) was known before the introduction of the monoid $\mathcal{P}$.  In 1993, in \cite{R} S.E. Rao provided a complete list of all simple objects in $\mathcal{F}^\sigma$, although the list is redundant -- no criterion is given to tell us when two simple objects are isomorphic.   But we provide here a sketch of this result, because the methods used there have been extended to the representation theory of multiloop algebras, which are natural generalizations of the loop algebras $L(\lie g)^\sigma$.  

Let $P \in \mathbb{C}\left[ t^{\pm m} \right]$ and $\left\langle P \right\rangle $ the ideal in $\mathbb{C}\left[ t^{\pm m} \right]$ generated by $P$.  Corresponding to any such $P$ is an ideal of $L(\lie g)^\sigma$ 
\[ PL(\lie g)^\sigma = \bigoplus_{\epsilon=0}^{m-1} \lie g_\epsilon \otimes t^\epsilon \left\langle P \right\rangle,\]
and we saw above in Section \ref{loopdef}, Proposition \ref{loopideal} that any ideal of $L(\lie g)^\sigma$ is of this form.  We denote by $L(\lie g)^\sigma(P)$ the quotient $L(\lie g)^\sigma / P L(\lie g)^\sigma$. Also recall that, for any $P \in \mathbb{C}\left[ t^m \right]$, $\dim(L(\lie g)^\sigma (P))< \infty$.

Let $\phi: L(\lie g)^\sigma \rightarrow \text{End}_{\mathbb{C}}(V)$ be a finite--dimensional irreducible representation of $L(\lie g)^\sigma$.  Then $\ker(\phi) \neq 0$, and therefore by Proposition \ref{loopideal} $\phi$ factors through some quotient 
\[ \xymatrix{
L(\lie g)^\sigma \ar[rr]^\phi \ar@{->>}[dd]& & \text{End}_{\mathbb{C}}(V)\\  \\ L(\lie g)^\sigma(P) \ar@{.>}_{\phi'}[uurr]& & }\]
where $P \in \mathbb{C}\left[ t^m \right]$.  By Lie's theorem, the representation $L(\lie g)^\sigma(P) \rightarrow \text{End}_{\mathbb{C}}(V)$ factors through the nilradical of $L(\lie g)^\sigma(P)$, which we denote $\nilrad(L(P))$.  Furthermore, the corresponding quotient of $L(\lie g)^\sigma(P)$ is semisimple:
\begin{prop} Let $P = \prod_{i=1}^\ell (t^m - a_i^m)^{n_i} \in \mathbb{C}\left[ t^m \right]$; $a_i^m \neq a_j^m$.  Then 
\[ L(\lie g)^\sigma(P) / \nilrad(L(P)) \cong \bigoplus_{i=1}^\ell \lie g,\]
and the quotient map $\pi: L(\lie g)^\sigma \twoheadrightarrow L(\lie g)^\sigma(P) / \nilrad(L(P))\cong \bigoplus_{i=1}^\ell \lie g$ is given by 
\[ \xymatrix@R=0.03in{\pi: L(\lie g)^\sigma \ar@{->>}[r] & \bigoplus_{i=1}^\ell \lie g\\
    x \otimes t^s \ar@{|->}[r]  & (a_1^s x , \ldots, a_\ell^s x).}\]
\end{prop}
Let us write $\lie g^\ell := \bigoplus_{i=1}^\ell \lie g$.  Therefore, for some $\ell \in \mathbb{Z}_{> 0}$ the homomorphism $\phi'$ from the preceding diagram factors through $\lie g^\ell$:
\begin{equation}\label{diagram2} \xymatrix{
L(\lie g)^\sigma \ar[rr]^\phi \ar@{->>}[dd]& & \text{End}_{\mathbb{C}}(V)\\  \\ L(\lie g)^\sigma(P) \ar@{->>}[rr] \ar@{.>}_{\phi'}[uurr]& & \lie g^\ell \ar@<2ex>@{.>}_{\phi''}[uu] }\end{equation}
and so $V$ is a finite--dimensional irreducible representation of $\lie g^\ell$.  Therefore as a $\lie g^\ell$--module, 
\[ V \cong_{\lie g^\ell} \otimes_{i=1}^\ell V(\lambda_i),\]
where $V(\lambda_i)$ is the finite--dimensional irreducible $\lie g$-module with highest weight $\lambda_i \in P^+$ (see {\cite[\S7, no.7]{B}}).  It is now evident that $V$ is an evaluation representation, as described in Section \ref{reps}.  

Recently, M. Lau has extended this result to the multiloop algebras in \cite{Lau}.  We now describe these algebras, and sketch the results contained there.  

Let $N \in \mathbb{Z}_{>0}$ and $R = \mathbb{C}\left[ t_1^\pm, \ldots, t_N^\pm\right]$.  The \textit{untwisted multiloop algebra} is the vector space $\lie g \otimes R$ with Lie bracket given by 
\[ \left[ x \otimes f, y \otimes g \right] = \left[x,y\right] \otimes fg,\]
for all $x, y \in \lie g$ and $f, g \in R$.  Let $\sigma_1, \ldots, \sigma_N \in \text{aut}(\lie g)$ be a collection of $N$ commuting automorphisms of finite order $m_i = \left| \sigma_i \right|$, and set $R_0 = \mathbb{C}\left[t_1^{\pm m_1}, \ldots, t_N^{\pm m_N} \right]$.  Let $G$ be the group 
\[ G = \mathbb{Z} / m_1\mathbb{Z} \oplus \cdots  \oplus \mathbb{Z} / m_N \mathbb{Z}   ,\]
and for each $i = 1, \ldots, N$ let $\zeta_i$ be a primitive $m_i^{th}$ root of unity.  Then $\lie g$ has a common eigenspace decomposition 
\[ \lie g = \bigoplus_{\overline{k} \in G} \lie g_{\overline{k}},\]
where $\overline{k}$ is the image of $k = (k_1, \ldots, k_N) \in \mathbb{Z}^N$ under the canonical map $\mathbb{Z}^N \twoheadrightarrow G$ and $\lie g_{\overline{k}} = \left\{ x \in \lie g : \sigma_i(x) = \zeta_i^{k_i}x \text{ for }i=1, \ldots, N \right\}$.  

The \textit{twisted multiloop algebra} $\mathcal{L} = L(\lie g, \sigma_1, \ldots, \sigma_N)$ is the Lie subalgebra 
\[ \mathcal{L} = \bigoplus_{k \in \mathbb{Z}^N} \lie g_{\overline{k}} \otimes \mathbb{C}t^k \subseteq \lie g \otimes R,\]
where $t^k := t_1^{k_1} \cdots t_N^{k_N}$.  The number of indeterminates $N$ in the Laurent polynomial ring $R$ is referred to in the literature as the \textit{nullity} of the multiloop algebra.  

There is a classification of the irreducible finite--dimensional representations of $\mathcal{L}$ which is a natural generalization of the $N=1$ case described above.  The first step toward this generalization was given by P. Batra in \cite{Bat}.  The multiloop algebras considered there are those for which $\sigma_1$ is a diagram automorphism, and $\sigma_2 = \cdots  = \sigma_N = \Id$. The case of an arbitrary multiloop algebra, addressed recently by M. Lau in \cite{Lau}, subsumes these previous results and is summarized here.

 We saw in the nullity 1 case that the action of $\mathcal{L}$ on $V$ factors through to a finite--dimensional semisimple Lie algebra.  The same is true in the more general nullity $N > 0$ case, but the result requires a new approach, since $R$ is no longer a PID for $N > 1$.  Therefore we cannot characterize the ideals of $\mathcal{L}$ via  a principally generated ideal $\left\langle P \right\rangle \subseteq \mathbb{C}\left[ t_1^{m_1}, \ldots, t_N^{m_N}\right]$, as was done in \cite{Kac} for $N=1$.  This phenomenon is replaced by the following generalization:

\begin{prop}[{\cite[Propositions 2.10, 2.13]{Lau}}]
There exists a radical ideal $I_0 \subseteq R_0$ such that 
\[\ker(\phi) = \bigoplus_{\overline{k} \in G} \lie g_{\overline{k}} \otimes t^k I_0.\]
The ideal $I_0$ has finite codimension, and there exists a unique set $\left\{a_1, \ldots, a_\ell \right\} \subseteq (\mathbb{C}^\times)^N$, $a_i = (a_{i1}, \ldots, a_{iN})$, such that $I_0 = \cap_{i=1}^\ell \lie{m}_{a_i}$, the intersection of the maximal ideals $\lie{m}_{a_i} = \left\langle (t_1 - a_{i1})\cdots (t_N - a_{iN}) \right\rangle \subseteq (\mathbb{C}^\times)^N$.  Furthermore, these points $a_i \in (\mathbb{C}^\times)^N$ satisfy $m(a_i) \neq m(a_j)$ for $1 \leq i \neq j \leq \ell$, where $m(a_i) = (a_{i1}^{m_1}, \ldots , a_{iN}^{m_N})$. 
\end{prop}

\begin{thm}[{\cite[Theorem 4.9]{Lau}}] Let $\phi: \mathcal{L} \rightarrow \text{End}_{\mathbb{C}}(V)$ be a finite--dimensional irreducible representation of $\mathcal{L}$.  Then
\[  \mathcal{L} / \ker(\phi) \cong \bigoplus_{i=1}^\ell \lie g\]
for some $\ell \in \mathbb{Z}_{>0}$. 
\end{thm}
\noindent The canonical projection $\mathcal{L} \twoheadrightarrow \mathcal{L} / \ker{\phi} \cong \lie g^\ell$ is given by an evaluation map
\begin{align*} \mathcal{L} & \twoheadrightarrow \lie g^\ell \\ x \otimes f & \mapsto (f(a_1)x, \ldots, f(a_\ell) x).\end{align*}
Hence the representation $\phi: \mathcal{L} \rightarrow \text{End}_{\mathbb{C}}(V)$ factors through $\lie g^\ell$ as above, and $V$ is again isomorphic (as a $\lie g^\ell$--module) to a tensor product $\otimes_{i=1}^\ell V(\lambda_i)$ of irreducible $\lie g$--modules.  Therefore $\phi$ is given by the composition 
\[ \xymatrix{ \mathcal{L} \ar@{->>}[r] & \lie g^\ell \ar[r] & \text{End}_{\mathbb{C}}(V(\lambda_1) \otimes \cdots \otimes V(\lambda_\ell)).}\]
Such a representation is sometimes denoted $V(\underline{\lambda}, \underline{a})$, where $\underline{\lambda} = ( \lambda_1, \ldots, \lambda_\ell ) \subseteq (P^+)^\ell$, and $\underline{a} = (a_1, \ldots, a_\ell)\subseteq ((\mathbb{C}^\times)^N)^\ell$. We summarize with one of the main theorems in \cite{Lau}: 

\begin{thm}[{\cite[Corollary 4.11, Theorem 4.12]{Lau}}]
Let $V$ be a finite--dimensional irreducible $\mathcal{L}$--module.  Then there exists an integer $\ell > 0$ and a pair $(\underline{\lambda}, \underline{a}) \in (P^+)^\ell \times ((\mathbb{C}^\times)^N)^\ell$, $\underline{a} = \left(a_1, \ldots, a_\ell\right)$, $m(a_i) \neq m(a_j)$ for $i \neq j$, such that 
\[ V \cong V(\underline{\lambda}, \underline{a}).\]
\end{thm}

And conversely to every pair $(\underline{\lambda}, \underline{a})$ satisfying these conditions corresponds a finite dimensional irreducible $\mathcal{L}$--module $V(\underline{\lambda}, \underline{a})$.   However, the collection of (isomorphism classes of) finite--dimensional irreducible $\mathcal{L}$--modules is not, in the general case,  parametrized by the collection of all such pairs $(\underline{\lambda}, \underline{a})$ (although this is true when $\sigma_i = \Id$ for all $i=1,\ldots, N$; i.e., when $\mathcal{L} = \lie g \otimes R$).  For example, let $\lie g$ be of type $A_3$, $N=1$ and $\sigma \in \aut(\lie g)$ the order 2 diagram automorphism of $\lie g$.  Set $\mathcal{L}= \lie g\otimes \mathbb{C}\left[t^\pm \right]$, and $\mathcal{L}^\sigma = L(\lie g, \sigma, 2)$.  Then for any $a \in \mathbb{C}^\times$, we have 
\begin{align*}
V(\omega_1, a) & \not \cong_{\mathcal{L}} V(\omega_3, -a), \; \; \text{ but }\\
V(\omega_1, a) &  \cong_{\mathcal{L}^\sigma} V(\omega_3, -a).\end{align*}
The final section in \cite{Lau} addresses this redundancy by providing an isomorphism criterion for the finite--dimensional irreducible $\mathcal{L}$--modules, telling us when and only when two such modules $V(\underline{\lambda}, \underline{a})$ and $V(\underline{\lambda}', \underline{a}')$ are isomorphic.

\section{Additional topics}\label{additional}

\subsection{Hyperloop algebra representation theory}\mbox{}

An integral form $\lie A_{\mathbb{Z}}$ of an algebra $\lie A$ over a field $\mathbb{F}$ is an algebra over $\mathbb{Z}$ such that $\lie A_{\mathbb{Z}}  \otimes_{\mathbb{Z}}\mathbb{F} = \lie A$.  In \cite{G}, H. Garland introduced an integral form of $U(L(\lie g))$.  This integral form allows for the construction of the \textit{hyperloop algebra} $U(L(\lie g))_{\mathbb{F}}$ defined over a field $\mathbb{F}$ of positive characteristic.  The finite--dimensional loop--highest weight representation theory of these algebras was initiated in \cite{JM} and \cite{JM2}, in which the simple objects are classified,  positive--characteristic analogs of the Weyl modules are constructed, and the block decomposition of the corresponding category is investigated.

\pagebreak 

\subsection{Infinite--dimensional representations}\mbox{}

The evaluation representations $V(\underline{\lambda}, \underline{a})$ were described in \cite{C, CPunitary} in 1986.  Here they appeared in the classification of certain infinite--dimensional representations of the Lie algebra $L(\lie g) \oplus \mathbb{C}d$, for which the bracket operation is defined by
\[ \left[ x \otimes t^r, y \otimes t^s \right] = \left[ x,y\right] \otimes t^{r+s}, \; \; \; \left[ d, x \otimes t^r \right] = r x \otimes t^s, \; \; \; \left[ d, d\right] = 0.\]
The element $d$ is referred to as `the' derivation of $L(\lie g)$, and $L(\lie g) \oplus \mathbb{C}d$ as the \textit{extended loop algebra}, denoted $L(\lie g)^e$.  The representations of these extended loop algebras are important in their own right - they are the so--called level zero representations of the affine Kac-Moody algebra $L(\lie g) \oplus \mathbb{C}c \oplus \mathbb{C}d$ for which the center $c$ acts trivially.  The representations considered in these references are those $L(\lie g)^e$--modules $V$ for which the $L(\lie g)^e$ action is locally finite and for which $\dim\left( V_{\lambda}\right) < \infty$ for all $\lambda \in (\lie h \oplus \mathbb{C}d)^*$ (i.e., the representations are integrable, with finite--dimensional weight spaces).  The category of all such representations is denoted $\mathcal{I}_{fin}$.

As a consequence of the derivation $d \in  L(\lie g)^e$, it can be shown that there are no non--trivial finite--dimensional representations of $L(\lie g)^e$.  In particular, $V(\underline{\lambda}, \underline{a})$ has no non--trivial $L(\lie g)^e$--module structure.  However, the vector space
\[ L(V(\underline{\lambda}, \underline{a})) = V(\underline{\lambda}, \underline{a}) \otimes \mathbb{C}\left[ t^{\pm 1} \right]\]
can be equipped with a non--trivial $L(\lie g)^e$--module structure (corresponding to this construction $V \mapsto L(V)$ is a natural functor $L: \mathcal{F} \mapsto \mathcal{I}_{fin}$; see \cite{CG}).  It was shown in \cite{C, CPunitary} that any such module $L(V(\underline{\lambda}, \underline{a}))$ is completely reducible, and that any simple object in $\mathcal{I}_{fin}$ occurs as a direct summand of such a module.  These results were extended to the twisted case in \cite{CPint}.  V. Chari and J. Greenstein showed in \cite{CG} that the blocks of the category $\mathcal{I}_{fin}$ are described by certain orbits of spectral characters under an action of $\mathbb{C}^\times$ on $\Xi$.  Irreducible integrable representations of the extended (untwisted) multiloops were classified in \cite{R2}, and recently this classification was extended to certain extended twisted multiloop algebras in \cite{PB}.  

\subsection{Finite--dimensional loop representation theory and equivariant functions}\mbox{}

  The loop algebra $L(\lie g)^\sigma$ may be identified with the Lie algebra of $\sigma$--equivariant regular functions $f: \mathbb{C}^\times \rightarrow \lie g$ (see {\cite[\S 8]{Kac}}).  This identification of the elements of $L(\lie g)^\sigma$ with functions on the complex torus extends to a parametrization of the simple objects and of the blocks of $\mathcal{F}^\sigma$; the details of these interpretations will be given in a forthcoming paper.  Also see {\cite[Corollary 5.20] {Lau}} for an application of this identification to the parametrization of irreducible representations of a multiloop algebra $\mathcal{L}$.  

\subsection{Current algebras and their finite--dimensional representations}\mbox{}

  Let $\lie g\left[ t \right] = \lie g \otimes \mathbb{C}\left[ t \right]$. This Lie subalgebra of $L(\lie g)$ is known as the \textit{current algebra} of $\lie g$, and there are a number of papers devoted to its finite--dimensional representations -- \cite{CG2, CL, CM2, FoL1, FoL}, for example.  In these papers, the Weyl modules are defined for $\lie g\left[ t \right]$, and have the same universal highest--weight property described as that described in this survey.  It is shown in \cite{FL2} that these Weyl modules for the current algebras are isomorphic to the Weyl modules for the loop algebras.  Other prominent examples of finite--dimensional current algebra representations found in these references include the Demazure and Kirillov--Reshetikhin modules (suitably defined as $\lie g\left[ t \right]$--modules).  All of these representations lie in the category of finite--dimensional $\mathbb{Z}_+$--graded $\lie g\left[ t\right]$--modules.  This category is studied in \cite{CG2}, where the authors show that many interesting representations of quivers arise from the study of such representations.

\subsection{Weyl modules of multiloop algebras}\mbox{}
The notion of a Weyl module can be generalized to any algebra for which the notion of a `highest--weight' representation is available.  Weyl modules for multiloop algebras and multivariable current algebras have been studied in \cite{FKL}, \cite{FL1}, \cite{FL2}, and \cite{Lok}.  There is a significant role played by these Weyl modules in the finite--dimensional representation theory of these algebras.  Analogs of Weyl modules for an arbitrary (twisted) multiloop algebra are so far unavailable, but would prove useful in the description of the corresponding category of finite--dimensional representations (via a block decomposition, as in \cite{CM}, \cite{S}).

\end{document}